\documentclass[12pt, reqno]{amsart}

\usepackage{graphicx}
\usepackage{mathrsfs}
\usepackage{color}
\usepackage{comment}
\usepackage{amssymb}
\usepackage{amsmath}
\usepackage{physics}
\usepackage{mathtools}
\usepackage{tikz, tikz-cd}
\usepackage
[colorlinks=true, urlcolor=blue, citecolor=blue, linkcolor=blue]{hyperref}
\usepackage{enumitem}

\usepackage{microtype} 
\usepackage{lmodern}
\usepackage{csquotes}

\calclayout

\let\oldtocsection=\tocsection
\let\oldtocsubsection=\tocsubsection
\let\oldtocsubsubsection=\tocsubsubsection

\renewcommand{\tocsection}[2]{\hspace{0em}\oldtocsection{#1}{#2}}
\renewcommand{\tocsubsection}[2]{\hspace{3em}\oldtocsubsection{#1}{#2}}
\renewcommand{\tocsubsubsection}[2]{\hspace{6em}\oldtocsubsubsection{#1}{#2}}

\numberwithin{equation}{section}

\newtheorem{prop}{Proposition}
\newtheorem{theorem}{Theorem}[]
\newtheorem*{theorem*}{Theorem}
\newtheorem{lemma}{Lemma}[section]
\newtheorem{corollary}[prop]{Corollary}

\newtheorem{theorema}{Theorem}

\theoremstyle{definition}
\newtheorem{definition}[]{Definition}

\newtheorem{quest}[prop]{Question}
\newtheorem{conjecture}[]{Conjecture}

\newtheorem{exam}[prop]{Example}
\newtheorem{remark}{Remark}

\newcommand{\dimmc}{\dim_{\mathrm{MC}}}
\newcommand{\mcov}{\widetilde{M}}
\newcommand{\zcov}{\widetilde{Z}}
\newcommand{\ez}{E_Z}
\newcommand{\ezcov}{\widetilde{E_Z}}
\newcommand{\lk}{L^{\otimes k}}

\newcommand{\slkz}{\mathbb{S}(L^{\otimes k}|_{Z})}

\newcommand{\diam}{\mathrm{diam}}

\newcommand{\QQ}{\mathbf{Q}}
\newcommand{\RR}{\mathbf{R}}
\newcommand{\ZZ}{\mathbf{Z}}

\begin{document}
\title{Circle bundles with PSC over large manifolds}



\author{Aditya Kumar}
\address{Department of Mathematics, University of Maryland, 4176 Campus Dr, College
Park, MD 20742, USA}
\email{akumar65@umd.edu}

\author{Balarka Sen}
\address{School of Mathematics, Tata Institute of Fundamental Research. 1, Homi Bhabha Road, Mumbai-400005, India}
\email{balarka2000@gmail.com, balarka@math.tifr.res.in}

\begin{abstract}
We construct infinitely many examples of macroscopically large manifolds of dimension $m \geq 4$ equipped with circle bundles whose total spaces admit metrics of positive scalar curvature and have macroscopic dimension at most $\lceil  m/2 \rceil + 1$. In particular, we answer a question of Gromov on the existence of circle bundles over enlargeable manifolds whose total spaces admit metrics of positive scalar curvature, in all dimensions. Our constructions are based on techniques from symplectic geometry. 
\end{abstract}
\maketitle

\tableofcontents

\section{Introduction} 
\subsection{Background and motivation}
\begin{definition}[Urysohn width]
    Let $X$ be a metric space. The $k$-Urysohn width of $X$, denoted by $\mathrm{UW}_k(X)$, is the infimum over $t > 0$ such that there is a $k$-dimensional simplicial complex $Y$ and a continuous map $f : X \to Y^k$ satisfying the property that the diameter of each fiber of $f$ is at most $t$, i.e. for all $y \in Y$, we have $\diam(f^{-1}\{y\} ) \leq t$.
\end{definition}

If $k<n$, then a small $k$-Urysohn width intuitively means that the manifold is \textit{close} to a $k$-dimensional space. This notion of closeness is captured in the following definition.

\begin{definition}[Macroscopic dimension]\cite{Gro96, drmoscow}\label{def-md}
    Let $X$ be a metric space. One says $\dim_{\mathrm{mc}}(X) \leq k$ (resp.~$\dim_{\mathrm{MC}}(X) \leq k$) if there is a continuous (resp.~Lipschitz) map $f: X \to Y^k$ to a $k$-dimensional simplicial complex $Y$ equipped with the uniform metric\footnote{The \emph{uniform metric} on a simplicial complex $K$ is one with respect to which all the $n$-simplices of $K$ are isometric to the standard Euclidean $n$-simplex $\Delta^n$, for all $n$.}, such that for all $y \in Y$, $\diam(f^{-1}(y)) < C$ for some uniform constant $C$.
\end{definition}

\begin{remark}
    Note that it is immediate that $\dim_{\mathrm{mc}}(X) \leq \dimmc(X)$\footnote{In general, these two notions may not agree even for universal cover of closed manifolds, see \cite{drmoscow}.}. For ease of exposition and simplicity, in this article we shall mostly concern ourselves with the slightly better behaved quantity $\dimmc(X)$ and refer to it as the \emph{macroscopic dimension} of $X$, except on occasions where we explicitly state otherwise. 

    Henceforth, by macroscopic dimension of a closed manifold $M$, we shall mean $\dimmc(\widetilde{M})$ where $\widetilde{M}$ is the universal cover of $M$ equipped with the lift of some metric on $M$. As any two metrics on a closed manifold are Lipschitz equivalent, this quantity is independent of the chosen metric on $M$. 
\end{remark}

The notions of Urysohn width and macroscopic dimension appear in several works of Gromov, with one of the aims being to capture the appropriate notion of smallness that must be exhibited by a manifold admitting a metric of positive scalar curvature (PSC). These also appear naturally in the study of non-collapsing and thick-thin decomposition. See also some of these recent works on the topic \cite{abg,guth,nabu,st}.
The primary motivation behind this article is the following conjecture formulated by Gromov for manifolds with positive scalar curvature \cite[Section $2 \text{\textonehalf}$]{Gro96}. 
\begin{conjecture}[Gromov] \label{conj}
    Let $M$ be a closed connected $n$-dimensional manifold, admitting a PSC metric. Then the macroscopic dimension of the universal cover of $M$ is at most $n-2$, that is, $$\mathrm{dim}_{\mathrm{mc}}(\widetilde{M}) \leq n-2 .$$
\end{conjecture}
Intuitively, the conjecture says that a manifold admitting a PSC metric, on large scales, looks like a manifold of at least two dimensions less. The primary evidence in this direction is the following simple observation: Given any closed manifold $M^{n-2}$, consider the manifold $M^{n-2} \times S^2_r$ with the product metric. Here, $S^2_r$ denotes a sphere of radius $r$. Then for $r$ small enough, the product metric on $M \times S^2_r$ has positive scalar curvature. In light of the previous remark, we pass to the universal cover of $M^{n-2} \times S^2_r$, which is $\widetilde{M} \times S^2_r$. Since the projection map to $\mcov$ is Lipschitz and each fiber is $S^2_r$, we have
$$\mathrm{UW}_{n-2}(\widetilde{M \times S^2_r}) \leq 2r,$$
and the macroscopic dimension $\dim_{\mathrm{mc}}(\widetilde{M \times S^2_r}) \leq \dimmc (\widetilde{M\times S^2_r}) \leq n-2$. Conjecture \ref{conj} asserts that all PSC manifolds must exhibit this behavior. 
Gromov's conjecture {was verified for $3$-manifolds in \cite{gl}. See also recent works \cite{cl25,lm23}} where a stronger version of Gromov's conjecture is established for $3$-manifolds. In higher dimensions there are several partial results under various assumptions (cf. \cite{bd, dd22}).  

While Conjecture \ref{conj} arises from the simple observation that existence of a PSC metric is guaranteed given the presence of an $S^2$ factor in the manifold, it is also natural to ask about the possible interaction between existence of PSC metrics and having an $S^1$ factor. This has been widely investigated. For instance, PSC manifolds with a trivial $S^1$ factor (i.e., $M \times S^1$) are amenable to the Schoen-Yau conformal descent argument \cite{sydescent}, as they have a non-trivial codimension-$1$ homology class given by the Poincar\'{e} dual of the $S^1$ factor. In fact, this was one of the key ideas in Schoen and Yau's proof that for $n\leq 7$, the torus  $T^n = S^1 \times \cdots \times S^1$ does not admit a PSC metric. Thereafter, Gromov and Lawson \cite{gl} proved this for all $n$. An important concept in their proof was the notion of enlargeability. We now recall the definition of enlargeable manifolds from \cite{gl}.

\begin{definition}[Enlargeability]
    A compact orientable Riemannian manifold $M$ of dimension $n$ is called \emph{$\varepsilon^{-1}$-hyperspherical} if there exists a $\varepsilon$-contracting map of positive degree from $M$ to the unit $n$-sphere. Further, $M$ is called \emph{enlargeable} if for every $\varepsilon > 0$, there exists a finite cover of $M$ {which} is $\varepsilon^{-1}$-hyperspherical and spin.
\end{definition}

The basic example of an enlargeable manifold is the unit circle $S^1 \cong \mathbf{R}/2\pi\mathbf{Z}$. Indeed, it is spin and the $n$-sheeted cover of $\mathbf{R}/2\pi\mathbf{Z}$ is $\mathbf{R}/2\pi n\mathbf{Z}$. The map $\mathbf{R}/2\pi n\mathbf{Z} \to \mathbf{R}/2\pi \mathbf{Z}$ given as scaling by $1/n$ is a $1/n$-contracting map, for any $n \geq 1$. It is straightforward to see that
\begin{enumerate}
    \item Product of enlargeable manifolds is enlargeable,
    \item A manifold $M$ admitting a positive degree {Lipschitz} map to an enlargeable manifold is also enlargeable, provided $M$ admits a finite spin covering.
\end{enumerate}
Therefore, $T^n \cong S^1 \times \cdots \times S^1$ is enlargeable and any spin manifold admitting a positive degree map to $T^n$ is also enlargeable. In \cite{gl}, it is proved that all compact hyperbolic manifolds are enlargeable. The notion of enlargeability is important because of the main result of Gromov-Lawson \cite[Theorem A]{gl}, which shows that enlargeable manifolds do not admit PSC metrics.

Note that if $M$ is enlargeable, then $M\times S^1$ is also enlargeable as the product of two enlargeable manifolds is enlargeable. However, the situation is less clear when the $S^1$ factor arises as a fiber of a \emph{non-trivial} circle bundle. In this direction, Gromov posed the following question \cite[pg. 661]{Gro18}:

\begin{quest}\cite[pg.~661]{Gro18}\label{qn-gro}
``What happens to nontrivial bundles $\underline{X} \to Y$? \newline
[...] In general, the examples in \cite{BH2009} indicate a possibility of non-enlargeable circle bundles over enlargeable $Y$; yet, it seems hard(er) to find such examples, where the corresponding $\underline{X}$ would admit complete metrics with $Sc \geq \sigma > 0$."
\end{quest}

In the negative direction we point out the recent comprehensive work due to He \cite{He25} on PSC metrics on non-trivial circle bundles. In particular, He shows that circle bundles with homologically non-trivial fibers do not admit PSC metrics if the base space is a $1$-enlargeable\footnote{$1$-enlargeable means that the maps to $S^n$ in the definition of enlargeability are of degree $1$.} manifold of dimension $n \leq 6$ and the circle fibers are homologically non-trivial in the total space \cite[Proposition 4.12]{He25}.

The main goal of this article is to study and exhibit examples of interactions between the existence of a positive scalar curvature metric and macroscopic dimension on non-trivial circle bundles over manifolds of dimension $n \geq 4$. In particular, we give a positive answer to Question \ref{qn-gro}. We construct these examples using techniques from symplectic geometry. The inspiration for our construction came from a remark due to Gromov \cite[pg.~8]{Gro23}:
\begin{displayquote}
    ``The emergent picture of spaces with $Sc.curv \geq 0$, where topology and geometry are intimately intertwined, is reminiscent of \emph{the symplectic geometry}, but the former has not reached yet the maturity of the latter."
\end{displayquote}
Granting this dictum of similarity between scalar curvature geometry and symplectic geometry, one is led to believe that the analogue of the Schoen-Yau descent for positive scalar curvature manifolds, in the symplectic case must be Donaldson's symplectic divisor theorem \cite{don96}, as it can be viewed as a codimension-$2$ descent-type result for symplectic manifolds. This analogy serves as the starting point of our constructions.

\subsection{Main results}
We discuss the necessary symplectic preliminaries in Section \ref{ddcon}.

For $n\geq 2$, in this article, we construct infinitely many examples of $2n$-manifolds $Z$, arising as Donaldson divisors (see, Theorem \ref{thm-dondiv}) in symplectic $(2n+2)$-manifolds, such that $Z$ does not admit a PSC metric and {macroscopic dimension of $Z$} is \textit{exactly} $2n$. Nevertheless, we construct circle bundles over these $Z$,
$$S^1 \hookrightarrow E_Z \xrightarrow{p} Z,$$
such that the $(2n+1)$-dimensional total space $E_Z$ admits a PSC metric. An intriguing feature of these examples is that $\dimmc(\ezcov) \leq n + 1 $.
Thus, even though $\ez$ is obtained from $Z$ by taking a twisted product with $S^1$, strangely the macroscopic dimension drops by half for $\ez$. For $n=2$, note that this is sharply in consonance with Gromov's macroscopic dimension conjecture (Conjecture \ref{conj}), i.e., $\dimmc(\ezcov) \leq \dim \ez-2$. We show this using the natural deformation retraction from $\ezcov$ to an $(n+1)$-dimensional simplicial complex arising as the spine of the complement of the Donaldson divisor $Z$ in the ambient symplectic $(2n+2)$-manifold. 

We remark that even though we state our results with $\dimmc$, similar techniques can be used to produce $2n$-manifolds $Z$ with $\dim_{\mathrm{mc}}(\zcov) = 2n$, which admit circle bundles whose total spaces exhibit an analogous macroscopic dimension drop phenomenon. We provide such examples in Example \ref{examplemdlarge}. Note that this phenomenon of the macroscopic dimension $\dim_\mathrm{mc}$ decreasing upon taking a (twisted) product with $S^1$ is not possible in the case of a trivial product. Indeed, \cite[Theorem 5.7]{drmoscow} shows that if  $\dim_\mathrm{mc}(\zcov) = 2n$, then 
$$\dim_\mathrm{mc}(\widetilde{Z \times S^1}) = 2n+1$$
To the best of our knowledge, examples exhibiting a drop in macroscopic dimension upon twisting with $S^1$ have not appeared before in literature. On a similar theme, we were informed by Guth \cite{guthemail} that a decrease in Urysohn width upon taking a double cover, i.e. an $S^0$ bundle, was observed by Balitskiy \cite[Theorem 2.4.4, pg.~36]{balthesis}.  

We now state the main result of this article:

\begin{theorema} \label{main}
    For $n\geq 2$, let $(M^{2n+2}, \omega)$ be a {closed} integral symplectic manifold. Let $Z^{2n}$ be {a} Donaldson divisor of $M$ with $[Z] = \mathrm{PD}(k [\omega])$, for some sufficiently large integer $k \gg 1$. Then there is a circle bundle 
    $S^1 \hookrightarrow E_Z \xrightarrow{p} Z,$ 
    such that the following statements hold:
    \begin{enumerate}[font=\normalfont]
        \item The total space $E_Z^{2n+1}$ admits a metric of positive scalar curvature. 
        \item The macroscopic dimension of $E_Z$ is at most $n+1$, i.e. $\dimmc(\ezcov) \leq n+1$.
        \item Let $(M^{2n+2}, \omega)$ be a symplectically aspherical\footnote{See Definition \ref{def-sympasph}. In particular, any symplectic manifold with $\pi_2 = 0$ (e.g. $T^{2n}$) is symplectically aspherical.} manifold with an amenable fundamental group. Then the macroscopic dimension of $Z^{2n}$ is exactly $2n$, i.e. $\dimmc(\widetilde{Z}) = 2n$.
        \item Let $M^{2n+2} = N^{2n} \times S_g$ and $\omega = \omega_N \oplus \omega_0$, where $\omega_0$ is an integral area form on $S_g$. If $N$ is enlargeable and $k$ is even, then $Z$ is also enlargeable. In particular, $Z$ does not admit a metric of positive scalar curvature.  
    \end{enumerate} 
    Additionally, varying the integer $k \gg 1$ produces an infinite collection of non-homeomorphic manifolds $Z^{2n}$ satisfying the above.
\end{theorema}

We emphasize that the Conclusions $(1)$ and $(2)$ {do} not rely on any assumptions either on the fundamental groups of $M$ and $\ez$, or on the (non-)existence of spin or almost spin structures on the universal cover of $\ez$. In particular, Conclusion $(2)$ does not follow from any existing results on Gromov's conjecture (Conjecture \ref{conj}). 

{Conclusion $(3)$ follows from combining Proposition \ref{prop-sympasphess} with a result of Dranishnikov \cite[Theorem 1.6]{dr11}. Conclusion $(4)$ follows from Proposition \ref{zpullm} and the fact that a manifold admitting a positive degree Lipschitz map to an enlargeable manifold is also enlargeable, provided it admits a finite cover that is spin.} Then, it follows from the classical result of Gromov and Lawson \cite{glsc} that $Z$ does not admit a PSC metric. Finally, it follows from Proposition \ref{prop-unique} that varying the parameter $k$ produces an infinite collection of such manifolds $Z$ satisfying the conclusions of Theorem \ref{main}.

Conclusion $(4)$ of Theorem \ref{main} gives the examples sought after in Gromov's question (Question \ref{qn-gro}) in all even dimension $\geq 4$. The case of odd dimensions $>4$ follows upon considering the bundle $p \times \mathrm{id} : E_Z \times S^1 \to Z \times S^1$ obtained by taking product of the bundle $p : E_Z \to Z$ with $S^1$, and observing that $Z\times S^1$ is enlargeable if $Z$ is enlargeable. 

Hanke \cite{hankemail} communicated to us a proof showing that total spaces of circle bundles over $3$-manifolds which are either enlargeable or macroscopic dimension $3$, do not admit positive scalar curvature metrics. We record the proof in Section \ref{subsec-thmb}.

\begin{theorema}\label{3dcase}
    Let $M$ be a closed connected $3$-manifold that is either enlargeable or has macroscopic dimension $3$. Let $p : E \to {M}$ be a circle bundle. Then the total space $E$ does not admit a metric of positive scalar curvature.
\end{theorema}

\begin{remark}
In dimension $2$, the only closed connected manifolds which are enlargeable or have macroscopic dimension $2$ are surfaces of genus $g \geq 1$, or their non-orientable quotients. Since these are aspherical, total spaces of circle bundles over such surfaces are also aspherical. By the classification of PSC $3$-manifolds, these can not admit positive scalar curvature metrics.

On the other hand, in dimension $1$, the total spaces of such circle bundles are either diffeomorphic to a torus or a Klein bottle. These can not admit a metric of positive scalar curvature either, by the Gauss-Bonnet theorem. 
\end{remark}

Theorem \ref{main} together with Theorem \ref{3dcase} provides a complete answer to Gromov's question (Question \ref{qn-gro}) for closed manifolds in all dimensions: if $\dim Y \geq 4$, there are infinitely many examples of circle bundles $\underline{X} \to Y$ where $Y$ is enlargeable and $\underline{X}$ admits a metric with $Sc > 0$. On the other hand, if $\dim Y \leq 3$, then there are no such examples.

\subsection{Organisation} 
In Section \ref{ddcon}, we define the necessary symplectic preliminaries and make some observations regarding Donaldson divisors. In Section \ref{sec-proof}, we prove Theorem \ref{main} and Theorem \ref{3dcase}. The reader may first read Section \ref{sec-proof}, and refer to Section \ref{ddcon} for the proofs of the results used in Section \ref{sec-proof} as they arise. In Section \ref{sec-eg}, we discuss several examples of manifolds satisfying the conclusions of Theorem \ref{main}, each illustrating different possible properties of the manifolds produces by our constructions.

\subsection*{Acknowledgements} The authors would like to thank Shihang He for several helpful comments on the first version of this paper, as well as for pointing out an error in a previous proof of Theorem \ref{3dcase}. They would also like to thank Alexey Balitskiy, Alexander Dranishnikov, and Larry Guth for answering their questions and providing several helpful remarks. The authors are also indebted to Bernhard Hanke for his invaluable input on an earlier version of the paper which helped simplify the proofs, as well as for communicating a proof of Theorem \ref{3dcase}. They also wish to thank the anonymous referees for pointing out inaccuracies as well as several suggestions that have improved this paper.  The second author wishes to thank his advisor Mahan Mj for his interest, as well as his friends and colleagues Ritwik Chakraborty and Sekh Kiran Ajij for several helpful conversations. The second author is supported by the Department of Atomic Energy, Government of India, under project no.12-R\&D-TFR-5.01-0500.

\section{Donaldson divisors: definition and properties}\label{ddcon} 

Let $(X, \omega)$ be a symplectic manifold. The symplectic form $\omega$ is said to be \emph{integral} if the cohomology class $[\omega] \in H^2(X; \RR)$ lies in the image of the change of coefficients homomorphism $H^2(X; \ZZ) \to H^2(X; \RR)$. 

\begin{remark}{
Any symplectic manifold $(X, \eta)$ admits a symplectic form $\omega$ which is integral. To see this, choose a collection of closed $2$-forms on $X$ representing a basis of the vector space $H^2(X; \RR)$, and express the class $\eta$ as an $\RR$-linear combination of the chosen closed $2$-forms. Since $H^2(X; \mathbf{Q}) \subset H^2(X; \mathbf{R})$ is dense, we may approximate the real coefficients by rational numbers and then multiply by an integer to clear out the denominators. Let $\omega$ be the resulting closed $2$-form. Then $[\omega] \in H^2(X; \mathbf{Z})$. Since non-degeneracy is both an open condition and is scale-invariant, $\omega$ is a non-degenerate $2$-form. Therefore, $\omega$ is an integral symplectic form, as desired.
}
\end{remark}

The following foundational result in symplectic geometry was established by Donaldson \cite{don96}. 

\begin{theorem}\cite[Theorem 1]{don96}\label{thm-dondiv}
Suppose $(X^{2n}, \omega)$ is a closed symplectic manifold where $\omega$ is integral. Then, for sufficiently large $k \gg 1$, there exists a symplectic submanifold $Z \subset X$ of codimension $2$, such that $Z$ represents the class Poincar\'{e} dual to $k[\omega]$.
\end{theorem}

The symplectic submanifolds $Z \subset (X, \omega)$ appearing in the statement of the theorem are referred to as \emph{Donaldson divisors} of $(X, \omega)$. We briefly outline the construction of $Z$. Since $\omega$ is an integral symplectic form, we can consider the \emph{pre-quantum} line bundle 
$$\pi : L \to X$$
for $(X, \omega)$. This is a complex line bundle with a hermitian connection $\nabla$ such that the curvature $2$-form is given by $F_{\nabla}=-2\pi i\omega$. Thus, note that the first Chern class of the line bundle $L$ is $c_1(L) = [\omega]$. The symplectic submanifold $Z^{2n-2} \subset (X^{2n}, \omega)$ is produced as the zero set of an \emph{asymptotically holomorphic} section of $L^{\otimes k}$. For details, see \cite[Theorem 5]{don96}. Note $c_1(L^{\otimes k})=k[\omega]$. Since $Z$ is the zero set of a section of $L^{\otimes k}$, we must have that $Z$ represents the class Poincar\'{e} dual to $k[\omega]$.

\subsection{Complement of Donaldson divisors}

We note the following two simple but useful observations.

\begin{lemma}\label{lem-nztriv}Let $\lk|_Z$ denote the restriction of $\lk$ to $Z$. Denote by $\slkz$ the unit circle bundle of $\lk|_Z$. Then the following statements are true. 
\begin{enumerate}[font=\normalfont]
    \item The tubular neighborhood $\nu(Z)$ of $Z \subset X$ is diffeomorphic to $L^{\otimes k}|_{Z}$ 
    \item The boundary $\partial \nu(Z)$ is diffeomorphic to $\slkz$.
\end{enumerate}
\end{lemma}

\begin{proof}
By the tubular neighborhood theorem, $\nu(Z)$ is diffeomorphic to the normal bundle of $Z \subset X$. By the construction of $Z$ outlined above, $Z = s^{-1}(\mathbf{0}),$ where $s : X \to L^{\otimes k}$ is a section transverse to the zero section $\mathbf{0} \subset L^{\otimes k}$. Since the normal bundle to $\mathbf{0} \subset L^{\otimes k}$ is $L^{\otimes k}$ itself, the normal bundle to $Z$ in $X$ is $L^{\otimes k}|_{Z}$ by transversality. Since $\nu(Z)\cong\lk_Z$, the total spaces of the unit circle bundles are also diffeomorphic, i.e. $\partial \nu(Z) \cong \slkz$.
\end{proof}

The following is an analogue of the Lefschetz hyperplane theorem in the symplectic setting, which we shall use in the next section:

\begin{prop}\cite[Proposition 39]{don96}\label{zpi12lemma}
The inclusion $j : Z^{2n-2} \hookrightarrow X^{2n}$ induces an isomorphism $j_* : \pi_k(Z) \to \pi_k(X),$ for all $k \leq n-2$, and a surjection for $k = n-1$.
\end{prop}

In fact, the following lemma can be deduced from the proof of Proposition \ref{zpi12lemma} in \cite[pg. 700-701]{don96}.

\begin{lemma}\label{lem-skel}
{The complement $X^{2n} \setminus Z$ admits a handlebody decomposition consisting of handles of index at most $n$. In particular, it deformation retracts to a simplicial complex of dimension $n$.}
\end{lemma}

\begin{proof}
By the construction of $Z \subset X$ summarized above, there is an asymptotically holomorphic section $s : X \to L^{\otimes k}$ such that $Z = s^{-1}( \mathbf{0})$. Let $\psi : X \setminus Z \to \RR$ be defined by $\psi := -\log |s|^2$. During the course of the proof of \cite[Proposition 39]{don96}, it is shown that $\psi$ is a function on $X \setminus Z$ with critical points having index at most $n$. {Therefore, by Morse theory, $X \setminus Z$ admits a handlebody decomposition consisting of handles of index at most $n$. Moreover, $X \setminus Z$ deformation retracts to the subcomplex of dimension at most $n$ given by union of the cores of the handles.}
\end{proof}

We deduce two corollaries from Lemma \ref{lem-skel}. These are central to the proof of Part (2) of Theorem \ref{main}.

\begin{corollary}\label{cormain}
    Let $\nu(Z)$ be a tubular neighborhood of $Z$ in $X^{2n}$. Let $W = X \setminus \nu(Z)$. Then the following statements are true,
    \begin{enumerate}[font=\normalfont]

        \item $W$ deformation retracts to a simplicial complex of dimension $n$.
        \item If $n \geq 3$, the inclusion $\partial W \hookrightarrow W$ induces an isomorphism on $\pi_1$.
    \end{enumerate}
\end{corollary}

\begin{proof}
For Part (1), observe $\partial W = \partial \nu(Z)$. So, by Lemma \ref{lem-nztriv}, we have $\partial W \cong \slkz$.   Therefore, $X \setminus Z$ is simply $W$ with a collar $\partial W \times [0, \infty)$ attached. Consequently, $X \setminus Z$ deformation retracts to $W$. Since by Lemma \ref{lem-skel}, $X \setminus Z$ deformation retracts to a $n$-dimensional simplicial complex, so does $W$. 

For Part $(2)$, {recall from Lemma \ref{lem-skel} that $W$ admits a handlebody decomposition consisting of handles of index at most $n$}. By reversing the handlebody decomposition, we get that $W$ can be constructed from $\partial W$ by attaching handles of index \emph{at least} $n$. If $n \geq 3$, then we conclude that all these handles are of index at least $3$. Since $\pi_1$ is unchanged under {attachments} of handles of index bigger than $2$, we conclude $\partial W \hookrightarrow W$ induces an isomorphism on $\pi_1$.
\end{proof}

\subsection{Topology of Donaldson divisors}

In the next proposition, we show that Donaldson divisors of a closed symplectic manifold $(X^{2n}, \omega)$ with integral symplectic form $\omega$, realizing the homology classes Poincar\'{e} dual to $k[\omega]$, are topologically distinct for distinct values of $k \gg 1$. We do this by calculating their Euler characteristic. 

\begin{prop}\label{prop-unique}
    Let $Z_k \subset (X^{2n}, \omega)$ denote a Donaldson divisor with $[Z_k] = \mathrm{PD}(k[\omega])$. Then, for sufficiently large values of $k \gg 1$, $|\chi(Z_k)| \sim O(k^n)$ is monotonically increasing in $k$. In particular, for $k, l \gg 1$ sufficiently large, $Z_k \cong Z_l$ if and only if $k = l$.
\end{prop}

\begin{proof}
    From Lemma \ref{lem-nztriv}, we know that the normal bundle of $Z_k$ in $X$ is isomorphic to $L^{\otimes k}|_{Z_k}$, where $L$ is a complex line bundle over $X$ with $c_1(L) = [\omega]$. Therefore,
    $$TZ_k \oplus L^{\otimes k}|_{Z_k} \cong TX|_{Z_k}$$
    Let $j : Z_k \hookrightarrow X$ denote the inclusion map. We compute the Chern classes:
    \begin{align}
        c_1(TZ_k) + c_1(L^{\otimes k}|_{Z_k}) &= j^*c_1(TX) \label{eq-1} \\
        c_i(TZ_k) + c_{i-1}(TZ_k) \smile c_1(L^{\otimes k}|_{Z_k}) &= j^*c_i(TX), \;\;\; i \geq 2 \label{eq-2}
    \end{align}
    Note that $c_1(L^{\otimes k}|_{Z_k}) = j^*c_1(L^{\otimes k}) = k \cdot j^*[\omega]$. Substituting this in (\ref{eq-1}), we obtain:
    \begin{equation}
        c_1(TZ_k) = j^* c_1(TX) - k \cdot j^*[\omega] \label{eq-3}
    \end{equation}
    Solving for the system of equations given by (\ref{eq-3}) and (\ref{eq-2}), we obtain
    \begin{equation}
        c_{n-1}(TZ_k) = \sum_{i = 0}^{n-1} {(-1)}^i k^i \cdot j^* c_{n-1-i} (TX) \smile [\omega]^{\smile i}, \label{eq-4}
    \end{equation}
where $c_0(TX) = 1$ by definition. Since the top Chern class is the Euler class, evaluating the cohomology classes appearing in either sides of (\ref{eq-4}) on $[Z_k] = k \cdot \mathrm{PD}[\omega]$ we get:
\begin{align*}
    \chi(Z_k) = \int_{Z_k} c_{n-1}(TZ_k) &= \int_X c_{n-1}(TZ_k) \smile k[\omega] \\ &= \sum_{i = 1}^{n} (-1)^{i-1} k^i \int_X {c_{n-i}(TX)} \smile [\omega]^{i}
\end{align*}
Note that $[\omega]$ and $c_i(TX)$, $0 \leq i \leq n$ are all independent of $k$. Hence, $\chi(Z_k)$ is a polynomial in $k$ with {degree $n$}, with coefficients independent of $k$, and the leading coefficient being nonzero. Therefore, for $k \gg 1$, $|\chi(Z_k)|$ is monotonically increasing in $k$, as required.
\end{proof}
\begin{remark}{The fact $|\chi(Z_k)| \sim O(k^n)$ should be put in context of the following standard observation: let $X = \mathbf{CP}^n$ and $V_k \subset \mathbf{CP}^n$ be a degree $k$ complex hypersurface. Then, 
$$\chi(V_k) = ((1-k)^{n+1}-1)k^{-1}+n+1.$$}
\end{remark}

\subsection{Donaldson divisors of \texorpdfstring{${N} \times S_g$}{N x Sg}}\label{ddms2}
Let $S_g$ be a closed surface with genus $g\geq 0$. Let $({N}^{2n}, \eta_0)$ be a closed integral symplectic manifold. We equip the manifold ${N} \times S_g$ with the integral symplectic form $\omega = \eta_0 \oplus \omega_0$, where $\omega_0$ denotes a positive, integral area form on $S_g$. We now make the following observation, which will be crucial to the proof of Part (4) of Theorem \ref{main}. 
\begin{prop} \label{zpullm}
The map $f : Z \to {N}$, obtained by the composition 
$$Z \hookrightarrow {N} \times S_g \xrightarrow{\pi} {N},$$
has positive degree.      
\end{prop}
\begin{proof}
    Take any regular value $p \in {N}$. Then the oriented count $\# f^{-1}(p)$ is equal to the oriented intersection number $\# ([Z] \cap [\{p\} \times S_g])$. Therefore, we have 
    \begin{align*}
        \deg(f) &= \# ([Z] \cap [\{p\} \times S_g]) \\
                            &= k \int_{\{p\} \times S_g}  \omega \quad \text{(since $Z$ is Poincar\'{e} dual to $k[\omega]$)} \\
                            &> 0 
    \end{align*}
    The last inequality holds because $\{p\} \times S_g$ is a symplectic submanifold of $({N}\times S_g, \omega)$. 
\end{proof}

\subsection{Symplectically aspherical manifolds}

\begin{definition}\label{def-sympasph}
    A symplectic manifold $(X, \omega)$ is said to be \emph{symplectically aspherical} if for every smooth map $\phi : S^2 \to X$, 
        $$\int_{S^2} \phi^*\omega = 0$$
\end{definition}

\begin{exam}\label{eg-sympash}
It is immediate from Definition \ref{def-sympasph} that symplectic submanifolds of symplectically aspherical manifolds are also symplectically aspherical. Since symplectic manifolds $(X, \omega)$ with $\pi_2(X) = 0$ (for instance, product of surfaces of genus $g \geq 1$) are automatically symplectically aspherical, Donaldson's theorem (Theorem \ref{thm-dondiv}) provides a large class of symplectically aspherical manifolds.
\end{exam}

This class of manifolds was originally introduced by Floer in the context of the Arnol'd conjecture in symplectic geometry. We refer the reader to the {survey} article \cite{KRTsurvey} for detailed information. The key property of symplectically aspherical manifolds that we shall be using is the following:

\begin{prop}\label{prop-sympasphess}
    Let $(X^{2n}, \omega)$ be a closed integral symplectic manifold which is symplectically aspherical. Then $X$ is a rationally essential manifold.\footnote{A $d$-dimensional manifold $M^d$ is rationally essential if for the map $\phi : M \to K(\pi_1(M), 1)$ classifying the universal cover of $M$, $\phi_*[M] \neq 0 \in H_d(K(\pi_1(M), 1); \QQ)$}
\end{prop}

\begin{proof}
    By attaching cells of dimension $\geq 3$ to $X$, we kill the higher homotopy groups $\pi_i(X)$, $i \geq 2$ to obtain a CW-complex which is a model for a $K(\pi_1(X), 1)$. By construction, there is a natural inclusion as a subcomplex
    $$i : X \hookrightarrow K(\pi_1(X), 1)$$
    Let $L$ be the pre-quantum line bundle over $X$, with $c_1(L) = [\omega]$. Since $(X, \omega)$ is symplectically aspherical, $L$ must extend to the $3$-skeleton of $K(\pi_1(X), 1)$. Indeed, for any $3$-cell of $K(\pi_1(X), 1)$, let $\phi : S^2 \to X$ denote the attaching map. Then, we have: 
    $$\langle c_1(\phi^* L), [S^2] \rangle = \int_{S^2} \phi^* \omega = 0$$
    Therefore, $c_1(\phi^* L) = 0$. Consequently, $\phi^* L$ must be the trivial bundle over the $2$-sphere. Thus, we may extend $L$ as a trivial bundle over the $3$-cell. Since complex line bundles over higher dimensional spheres are always trivial, we can extend $L$ over all the other higher dimension cells as well. 

    In this way, we obtain a line bundle $L'$ over $K(\pi_1(X), 1)$ such that $i^* L' = L$. Let us denote $\alpha := c_1(L') \in H^2(K(\pi_1, 1); \mathbf{Z})$. Then, $i^* \alpha = [\omega]$. Consequently, $i^*{\alpha^{\smile n}} = [\omega^{\wedge n}] \neq 0$. Thus, $i^* : H^{2n}(K(\pi_1(X), 1); \mathbf{Z}) \to H^{2n}(X; \mathbf{Z})$ is nonzero. Since $[\omega^{\wedge n}]$ is not torsion, ${\alpha^{\smile n}}$ must not be torsion either. Therefore, $i_* : H^{2n}(K(\pi_1, 1); \QQ) \to H^{2n}(X; \QQ)$ is nonzero. Therefore, by the universal coefficients theorem, $i_* : H_{2n}(X; \QQ) \to H_{2n}(K(\pi_1(X), 1); \QQ)$ must be nonzero as well.
\end{proof}

\section{Proofs of the main results}\label{sec-proof}

\subsection{Proof of Theorem \ref{main}}\label{subsec-thma}
In this section, we shall prove our main result. We begin with a general lemma relating macroscopic dimension of a manifold with the homotopy type of a nullcobordism of the manifold.

\begin{lemma} \label{lem-ret}
    Let $W^{{m}}$ be a compact ${m}$-dimensional manifold with boundary $E = \partial W$. Suppose
    \begin{enumerate}[font=\normalfont]
        \item $W$ deformation retracts to a $k$-dimensional simplicial complex $K \subset W$,
        \item The inclusion $i : E = \partial W \hookrightarrow W$ induces an isomorphism on $\pi_1$.
    \end{enumerate}
    Then, $\dimmc(\widetilde{E\phantom{.}}\!) \leq k$.
\end{lemma}

\begin{proof}
Let $r : W \to K$ denote the deformation retract. Let $f = r \circ i : E \to K$ denote the composition with the inclusion $i : E \hookrightarrow W$. Since $i$ and $r$ both induce an isomorphism on $\pi_1$, so does $f$. {Let $q_E : \widetilde{E\phantom{.}}\! \to E$ and $q_K : \widetilde{K} \to K$ denote the universal covers of $E$ and $K$ respectively. Since $f$ induces an isomorphism on the fundamental group, it lifts to a map between the respective universal covers:}
    \begin{gather*}
    \widetilde{f} : \widetilde{E\phantom{.}}\! \to \widetilde{K},
    \end{gather*}
{such that $q_K \circ \widetilde{f} = f \circ q_E$.} Choose an arbitrary Riemannian metric $g$ on ${E}$, and let $\widetilde{g}$ be its lift to the universal cover ${\widetilde{E\phantom{.}}}\!$. We will show that the fibers of $\widetilde{f}$ are uniformly bounded.

By the simplicial approximation theorem, we may assume without loss of generality that $f$ is a simplicial map with respect to some triangulation of $E$ and a subdivision of $K$. Note that $K$ is compact as it is a retract of a compact manifold $W$. Therefore, $K$ contains finitely many vertices $v_\alpha \in K^{(0)}$, $\alpha \in I$ where $I$ is a finite set of indices. Let $B_\alpha \subset K$ denote the open star of $v_\alpha$. By passing to an even finer subdivision of $K$ if necessary, we may assume without loss of generality that $B_\alpha$ is contractible for all $\alpha \in I$. Let $U_\alpha := f^{-1}(B_\alpha) \subset E$. 

We claim that the universal cover $q_E : \widetilde{E\phantom{.}}\! \to E$ trivializes over $U_\alpha$ for all $\alpha \in I$. It suffices to show that the inclusion $U_\alpha \hookrightarrow E$ induces the zero map on the fundamental group. To this end, let $\gamma \subset U_\alpha$ be a loop. Then $f_*[\gamma] = 0 \in \pi_1(B_\alpha)$ since $B_\alpha$ is contractible. As $f$ induces an isomorphism on $\pi_1$, $[\gamma] = 0 \in \pi_1(E)$. This proves the claim.

Let $\widetilde{K}$ be equipped with the simplicial structure obtained from lifting the simplicial structure on $K$. For each $\alpha \in I$, we denote the vertices of $\widetilde{K}$ contained in the fiber $q_K^{-1}(v_\alpha)$ as $w_{\alpha,\beta}$, indexed by $\beta \in J$. Then 
$$q_E^{-1}(U_\alpha) = \bigsqcup_{\beta\in J} V_{\alpha,\beta},$$
where $f(V_{\alpha,\beta}) \subset \widetilde{K}$ contains $w_{\alpha,\beta}$ and the restriction $q_E : V_{\alpha,\beta} \to U_\alpha$ is a diffeomorphism. In particular, $q_E : (V_{\alpha,\beta}, \widetilde{g}) \to (U_\alpha, g)$ is a Riemannian isometry for all $\alpha \in I, \beta \in J$. Therefore, the intrinsic diameter of $V_{\alpha, \beta}$ is equal to the intrinsic diameter of $U_\alpha$ for all $\alpha \in I, \beta \in J$. Let $C > 0$ denote the maximum of the intrinsic diameters of $(U_\alpha, g)$, $\alpha \in I$. This exists, as $I$ is a finite set. Since the extrinsic diameter is at most the intrinsic diameter for any open subset in a Riemannian manifold, we obtain $\diam(V_{\alpha, \beta}) \leq C$ for all $\alpha \in I, \beta \in J$.

Let $x \in \widetilde{K}$. Since $\{B_\alpha\}$ is an open cover of $K$, $q_K(x) \in B_\alpha$ for some $\alpha \in I$. Therefore, 
$$\widetilde{f}^{-1}(x) \subset \widetilde{f}^{-1}(q_K^{-1}(B_\alpha)) = q_E^{-1}(f^{-1}(B_\alpha)) = q_E^{-1}(U_\alpha)$$
Thus, there exists $\beta \in J$ such that $\widetilde{f}^{-1}(x) \subset V_{\alpha,\beta}$. Hence, $\diam(\widetilde{f}^{-1}(x)) \leq \diam(V_{\alpha,\beta}) \leq C$. Note that $C$ only depends on the simplicial structure on $K$, the map $f$ and the Riemannian manifold $(E, g)$. Thus, the fibers of $\widetilde{f}$ are uniformly bounded. {Furthermore, since $f$ has compact domain and target, and is simplicial, it is $\lambda$-Lipschitz for some $\lambda > 0$. Thus, $\widetilde{f}$ is also $\lambda$-Lipschitz.} Therefore, $\dimmc (\widetilde{E\phantom{.}}\!) \leq k$.\end{proof}

We are now ready to give a proof of Theorem \ref{main}. 
\begin{proof}[Proof of Theorem \ref{main}]
Let $(M^{2n+2}, \omega)$ be a closed integral symplectic manifold. Let $L$ be the pre-quantum line bundle on $(M , \omega)$. Let $Z^{2n} \subset (M , \omega)$ be a Donaldson divisor, with $[Z] = \mathrm{PD}(k \cdot [\omega])$, for $k \gg 1$ sufficiently large. We set $\ez := \slkz$ as the unit circle bundle of $L^{\otimes k}$ restricted to $Z$ (as introduced in Lemma \ref{lem-nztriv}). 

\begin{enumerate}[topsep=1ex,itemsep=1ex,partopsep=1ex,parsep=1ex, label=Proof of (\arabic*)., wide=0pt, font=\itshape]

\item Let $\nu(Z)$ be a tubular neighborhood of $Z$ in $X$. Consider the manifold with boundary $W := M \setminus \nu(Z)$. From Lemma \ref{lem-skel}, $W$ has a handle decomposition consisting of handles of {index $\leq n+1$}. Let $V := W \setminus {B^{2n+2}}$ be the $2n$-manifold obtained from removing a {$(2n+2)$}-ball $B^{2n+2} \subset W^\circ$ from the interior of $W$. Then $V$ is a cobordism between {$S^{2n+1}$} and $\partial \nu(Z)$. But ${\partial \nu(Z)= \mathbb{S}(L^{\otimes k}|_Z)}$ by Lemma \ref{lem-nztriv}. Note that $V$ has handles of index at most {$n+1$}. Therefore, $\slkz$ is obtained from ${S^{2n+1}}$ by surgeries of codimension $\geq {n+1} \geq 3$. Since ${S^{2n+1}}$ is PSC, so is $\slkz$ by the Gromov-Lawson-Schoen-Yau surgery theorem \cite{sydescent,glsc} which states that any manifold obtained from a PSC manifold by surgeries of {codimension} $\geq 3$ is also PSC.
\item Once again, let us consider the manifold with boundary $W := M \setminus \nu(Z)$ where $\nu(Z)$ denotes a tubular neighborhood of $Z$ in $X$. Note that $\partial W = \partial \nu(Z)$. By Lemma \ref{lem-nztriv}, $\partial \nu(Z) \cong \slkz = \ez$. Therefore $\partial W \cong \ez$. In Corollary \ref{cormain}, we deduced that $(W^{2n+2}, \ez)$ satisfies the hypothesis of Lemma \ref{lem-ret} with {$m=2n+2$ and $k = n+1$}. Therefore, we conclude $\dimmc(\ezcov) \leq n+1$.
\item Since $(M, \omega)$ is symplectically aspherical and $Z \subset M$ is a symplectic submanifold, $Z$ is also symplectically aspherical (see, Example \ref{eg-sympash}). Therefore, by Proposition \ref{prop-sympasphess}, $Z$ is a rationally essential manifold. By {Proposition \ref{zpi12lemma}}, $\pi_1(Z) \cong \pi_1(M)$. By hypothesis, $\pi_1(M)$ is amenable. Therefore, $Z^{2n}$ is a rationally essential manifold with amenable fundamental group. {Hence, by \cite[Theorem 1.6]{dr11}, $\dimmc(Z) = 2n$.} 
\item Here $M = N\times S_g$, with symplectic form $\omega = \omega_N \oplus \omega_0$, for some closed integral symplectic manifold $(N, \omega_N)$ and a positive integral area form $\omega_0$ on $S_g$. Choose a Donaldson divisor $Z \subset (N \times \Sigma_g, \omega)$ such that $[Z] = \mathrm{PD}(k \cdot [\omega])$, for some $k \gg 1$ sufficiently large such that, moreover, $k$ is even. Let $f : Z \to N$ be the restriction of the projection map $M = N \times S_g \to N$ to $Z$. By Proposition \ref{zpullm} the map $f : Z \to N$ is of positive degree. {By the Whitney approximation theorem, we may homotope $f$ slightly to be a smooth map without altering the degree. Since $Z$ is compact, the resulting map is $C$-Lipschitz for some $C > 0$.} Therefore, to show that $Z$ is also enlargeable, we just need to show that there is a finite cover $\widehat{Z} \to Z$ such that $\widehat{Z}$ is spin. 

We begin by noting that the inclusion $j : Z \hookrightarrow M$ induces an isomorphism on $\pi_1$, by Proposition \ref{zpi12lemma}. Since $N$ is enlargeable, there is a finite cover $\widehat{N} \to N$ such that $\widehat{N}$ is spin. Let $\widehat{M} := \widehat{N} \times S_g$. Since $\widehat{M}$ is a product of spin manifolds, it is also spin. As $p : \widehat{M} \to M$ is a finite cover, it corresponds to some finite index subgroup $H \leq \pi_1(M) \cong \pi_1(Z)$. Let $p : \widehat{Z} \to Z$ be the finite cover of $Z$ corresponding to $H \leq \pi_1(Z)$. Note that the inclusion $j : Z \hookrightarrow M$ lifts canonically to an inclusion $\widehat{j} : \widehat{Z} \hookrightarrow \widehat{M}$. 

Let $\widehat{\omega} := p^*\omega$ denote the lift of the symplectic form on $M$ to $\widehat{M}$. Then, $(\widehat{M}, \widehat{\omega})$ is a symplectic manifold, and $\widehat{j} : \widehat{Z} \to \widehat{M}$ is a symplectic embedding. From (\ref{eq-3}) in Proposition \ref{prop-unique}, we have:
\begin{equation}c_1(Z) = j^* c_1(M) - k \cdot j^*[\omega]\label{eq-chern}\end{equation}
By pulling back (\ref{eq-chern}) to the cover $p : \widehat{M} \to M$, we obtain:
\begin{equation}c_1(\widehat{Z}) = {\widehat{j}^* c_1(\widehat{M})} - k \cdot \widehat{j}^*[\widehat{\omega}]\label{eq-cherncov}\end{equation}
Next, we reduce both sides of (\ref{eq-cherncov}) modulo $2$. Recall that the first Chern class of an almost complex manifold is an integral lift of the second Stiefel-Whitney class. Since $\widehat{M}$ is spin, $w_2(\widehat{M}) = 0$. Furthermore, as $k$ is even, we conclude
$$w_2(\widehat{Z}) \equiv c_1(\widehat{Z}) \equiv 0 \pmod{2}$$
Hence, $\widehat{Z}$ is spin. Therefore, $Z$ is enlargeable. In particular, by Gromov-Lawson \cite{glsc}, $Z$ does not admit a metric of positive scalar curvature.
\end{enumerate}

Finally, it follows from Proposition \ref{prop-unique} that for distinct, {sufficiently} large integers $k \neq l$, the corresponding Donaldson divisors $Z_k, Z_l \subset (M, \omega)$ are not homeomorphic. Thus, one finds infinitely many manifolds $Z^{2n}$ satisfying Conclusions $(1)$-$(4)$. \end{proof}

\subsection{Proof of Theorem \ref{3dcase}}\label{subsec-thmb} We begin by stating the following proposition. The proof of this result was communicated to us by Hanke \cite{hankemail}.

\begin{prop}\label{prop-asphsum}
    Let $M^3$ be a closed orientable $3$-manifold such that $M$ is either enlargeable or $\dimmc(\widetilde{M}) = 3$. Then the Kneser-Milnor prime decomposition of $M$ must contain an aspherical summand.
\end{prop}

\begin{proof}
We prove the contrapositive. Let $M^3$ be a closed $3$-manifold. Let 
$$M = P_1 \# P_2 \# \cdots \# P_n$$
be the prime decomposition of $M$, and suppose $P_i$ is not aspherical for any $1 \leq i \leq n$. It follows from the sphere theorem that every closed orientable irreducible $3$-manifold with infinite fundamental group is aspherical. Therefore, either $P_i$ is prime but not irreducible so that $P_i \cong S^1 \times S^2$, or $P_i$ has finite fundamental group. Consequently, we must have $\pi_1(M) \cong F_k  * G_1 * \cdots * G_\ell$ where $G_i$ is finite for all $1 \leq i \leq \ell$ and $F_k$ is the free group of rank $k$. Therefore,
$$K(\pi_1(M), 1) \cong \bigvee_{i = 1}^k S^1 \vee \bigvee_{i 
= 1}^\ell K(G_i, 1)$$
Since $G_i$ is a finite group, $\widetilde{H}_*(K(G_i, 1); \mathbf{Q}) = 0$. Thus, $H_3(K(\pi_1(M), 1); \mathbf{Q}) = 0$. Therefore, $M$ can not be rationally essential. By \cite[Theorem 1.4]{hs06}, closed enlargeable manifolds are rationally essential. Therefore, $M$ can not be enlargeable either.

Next, we show $\dimmc(\widetilde{M}) \neq 3$. Consider the homomorphism
$$\varphi : \pi_1(M) \cong F_k * G_1 * \cdots * G_\ell \to G_1 \times \cdots \times G_\ell,$$
given by quotienting out the free factor $F_k$, followed by abelianization. Since $\varphi$ is surjective and $G_1 \times \cdots \times G_\ell$ is finite, $\ker(\varphi)$ must be a finite index subgroup of $\pi_1(M)$. It follows from the Kurosh subgroup theorem that $\ker(\varphi)$ is free of finite rank. Let $\widehat{M} \to M$ be the finite sheeted cover corresponding to $\ker(\varphi)$. 

Since $\pi_1(\widehat{M}) \cong \ker(\varphi)$ is free, any prime summand of $\widehat{M}$ is a prime $3$-manifold with free fundamental group. Note that there are no closed orientable irreducible $3$-manifold with free fundamental group. Indeed, such a $3$-manifold must be aspherical by the sphere theorem. But this implies the cohomological dimension of a free group is $3$, which is absurd. Therefore, all nontrivial prime summands of $\widehat{M}$ must be prime but not irreducible $3$-manifolds, i.e. $S^1 \times S^2$. Therefore, $\widehat{M} \cong \#^m (S^1 \times S^2)$. 

Next, observe that the universal cover of $M$ agrees with the universal cover of $\widehat{M}$. However, the universal cover of $\widehat{M} \cong \#^m (S^1 \times S^2)$ is obtained from the standard Cayley graph of $F_m$ by assigning vertex spaces as $S^3$, edge spaces as $S^2 \times I$, and tubing the two vertex spaces corresponding to the endpoints of any given edge by an $S^2 \times I$. Consider the map $f : \widetilde{M} \to Cay(F_k)$ obtained by collapsing down the vertex spaces to the vertices and edge spaces to the edges of the Cayley graph. Thus, $f$ is Lipschitz and the fibers of $f$ have uniformly bounded diameter. Therefore, $\dimmc(\widetilde{M}) \leq 1$. In particular, $\dimmc(\widetilde{M}) \neq 3$, as desired.
\end{proof}

We record a lemma that will be used in the proof of Theorem \ref{3dcase}. The proof of the lemma follows directly from the proof of \cite[Theorem 2.1]{bol06}. We provide the details below for completeness. For background on characteristic classes, we refer the reader to \cite{milstash}. See also the proof of \cite[Theorem 2.1]{bol06} which highlights the key notions and properties used. We adopt the following convention: given a class $\alpha \in H^i(X)$, we denote $|\alpha| = i$. 

\begin{lemma}\label{lem-nullfibers}
    Let $M$ be an orientable $3$-manifold such that its Kneser-Milnor prime decomposition contains an aspherical summand. Let $p : E \to M$ be a fiberwise orientable circle bundle. Suppose the fibers of $p$ are $\pi_1$-null in the total space of $E$. Then the classifying map $f_E : E \to K(\pi_1(E), 1)$ for the universal cover of $E$, can not be homotoped into the $2$-skeleton of $K(\pi_1(E), 1)$.
\end{lemma}

\begin{proof}
    Applying the homotopy fiber long exact sequence to $S^1 \hookrightarrow E \stackrel{p}{\to} M$, we obtain
    \begin{equation*}
    \pi_2(M) \to \pi_1(S^1) \to \pi_1(E) \stackrel{p_*}{\to} \pi_1(M) \to 0
    \label{eq-les}
    \end{equation*}
    Since $\pi_1(S^1) \to \pi_1(E)$ is zero by hypothesis, $p$ must induce an isomorphism on $\pi_1$. Let $q : D(E) \to M$ denote the $D^2$-bundle over $M$ associated to the circle bundle $p : E \to M$, so that $E = \partial D(E)$ and $q|_E = p$. Since $q$ is a homotopy equivalence, the inclusion $E = \partial D(E) \hookrightarrow D(E)$ induces an isomorphism on $\pi_1$. Thus, let us denote
    $$G := \pi_1(E) \cong \pi_1(D(E)) \cong \pi_1(M)$$
    Let $f_M : M \to K(G, 1)$ and $f_{D(E)} : D(E) \to K(G, 1)$ denote the classifying maps for the universal covers of $M$ and $D(E)$, respectively. Then $f_{D(E)}$ restricted to $E = \partial D(E) \subset D(E)$ is homotopic to $f_E$, and $f_{D(E)}$ restricted to the zero section $M \subset D(E)$ is homotopic to $f_M$.
    
    From the hypothesis, it follows that $M = N \# K$ for an orientable aspherical $3$-manifold $K$. Let $\varphi : K(G, 1) \to K(\pi_1(K), 1) \simeq K$ be induced from the quotient map to the free factor $G \cong \pi_1(N) * \pi_1(K) \to \pi_1(K)$. Let $\psi : K \to S^3$ be a degree $1$ map, given by collapsing the complement of a ball in $K$. Then, let $g : K(G, 1) \to S^3$ by defined by $g := \psi \circ \varphi$.
    
    If $f_E$ was homotopic into the $2$-skeleton of $K(G, 1)$, then $g \circ f_E$ would be null-homotopic. Therefore, $g \circ f_{D(E)}$ restricted to $E = \partial D(E) \subset D(E)$ would be null-homotopic. Hence, $g \circ f_{D(E)}$ factors through the quotient $D(E) \to D(E)/\partial D(E)$ up to homotopy. Note that $T := D(E)/\partial D(E)$ is the Thom space of the disk bundle $D(E)$. Thus, we obtain a map $F : T \to S^3$. Composing the zero section $M \to D(E)$ with the quotient map $D(E) \to T$, we obtain an embedding $s : M \to T$. We observe $F \circ s = g \circ f_{D(E)} \circ s$ is homotopic to $g \circ f_M$. Thus, $F \circ s : M \to S^3$ is degree $1$, as well. 

    Henceforth, we closely imitate the proof of \cite[Theorem 2.1]{bol06}. Let 
    $$\Phi : H^k(M; \mathbf{Z}_2) \to H^{k+2}(D(E), \partial D(E); \mathbf{Z}_2) \cong H^{k+2}(T; \mathbf{Z}_2)$$
    denote the Thom isomorphism, defined by $\Phi(x) = q^*x \smile u$ where $u \in H^2(T; \mathbf{Z}_2)$ denotes the Thom class of the disk bundle $q : D(E) \to M$. Let $\eta \in H^3(S^3; \mathbf{Z}_2) \cong \mathbf{Z}_2$ be the generator and let $a = F^*(\eta) \in H^3(T; \mathbf{Z}_2)$. Since $\Phi$ is an isomorphism, there exists $z \in H^1(M; \mathbf{Z}_2)$ such that $\Phi(z) = a$. Let $w_2 \in H^2(M; \mathbf{Z}_2)$ denote the second Stiefel-Whitney class of $q: D(E) \to M$. Using $s^* u = w_2$ and the fact that $F \circ s : M \to S^3$ is degree $1$, we obtain
    \begin{equation}
    z \smile w_2 = s^*q^* z\smile s^*u = s^*(q^*(z) \smile u) = s^*(\Phi(z)) = s^*(a) = (F \circ s)^*(\eta) \neq 0
    \label{eq-sq0}
    \end{equation}
    Let $\mathrm{Sq}^i$ denote the Steenrod square operations. Since $|z| = 1$, we get $\mathrm{Sq}^1(z) = \mathrm{Sq}^2(z) = 0$. Using Cartan's formula, naturality of Steenrod squares and $|u| = 2$, we compute
    \begin{equation}
    \mathrm{Sq}^2(a) = \mathrm{Sq}^2(\Phi(z)) = \mathrm{Sq}^2(q^*z \smile u) = q^*z \smile \mathrm{Sq}^2(u) = q^* z \smile u \smile u
    \label{eq-sq1}
    \end{equation}
    From Thom-Wu formula, we have $\Phi(w_2) = u \smile u$. Plugging this back in (\ref{eq-sq1}), we obtain
    \begin{equation}\mathrm{Sq}^2(a) = q^*z \smile \Phi(w_2) = q^*z \smile q^* w_2 \smile u = q^*(z \smile w_2) \smile u = \Phi(z \smile w_2)
    \label{eq-sq2}
    \end{equation}
    From (\ref{eq-sq0}), $z \smile w_2 \neq 0$. Since $\Phi$ is an isomorphism, combining this with (\ref{eq-sq2}) implies $\mathrm{Sq}^2(a) \neq 0$. On the other hand, by naturality of Steenrod squares, we have
    $$\mathrm{Sq}^2(a) = \mathrm{Sq}^2(F^* \eta) = F^* \mathrm{Sq}^2(\eta) = 0,$$
    since $\mathrm{Sq}^2(\eta) \in H^5(S^3; \mathbf{Z}_2)$ vanishes as $H^5(S^3; \mathbf{Z}_2) \cong 0$. This gives us the required contradiction. Therefore, $f_E$ is not homotopic into the $2$-skeleton of $K(G, 1)$, as desired. 
\end{proof}

We now proceed to the proof of Theorem \ref{3dcase}. As a crucial ingredient, we shall use a result due to He \cite[Theorem 1.9]{He25}.

\begin{proof}[Proof of Theorem \ref{3dcase}]
    Let $M^3$ be a closed $3$-manifold such that $M$ is either enlargeable or $\dimmc(\widetilde{M}) = 3$, and let $p : E \to M$ be a circle bundle over $M$ such that the total space of $E$ admits a PSC metric. As a first step, we argue that it can be ensured without loss of generality that $(1)$ $M$ is orientable, $(2)$ the circle bundle $E \to M$ is fiberwise orientable.
    
    \begin{enumerate}
        \item If $M$ is enlargeable, then $M$ is orientable by definition. If $\dimmc(\widetilde{M}) = 3$ and $M$ is not orientable, we pass to the orientation double cover of $M$ and pull the bundle $E \to M$ back to this cover. As the macroscopic dimension of a manifold is a quantity that only depends on the universal cover of the manifold, the cover continues to have macroscopic dimension $3$. Moreover, the total space of the pullback of $E \to M$ is a double cover of the total space of $E$. Therefore, the total space of the pulled back bundle also admits a PSC metric. Therefore, without loss of generality, we may assume $M$ is orientable. 
        
        \item If the circle bundle $E$ is not fiberwise orientable, then $w_1(E) \neq 0 \in H^1(M; \mathbf{Z}_2)$. By using the isomorphism $H^1(M; \mathbf{Z}_2) \cong \mathrm{Hom}(\pi_1(M), \mathbf{Z}_2),$
        we obtain a homomorphism $\phi : \pi_1(M) \to \mathbf{Z}_2$ corresponding to $w_1(E)$. We pass to the double cover of $M$ corresponding to the index $2$ subgroup $\ker(\phi)$ of $\pi_1(M)$. The pullback of $E$ to this cover is an oriented circle bundle by construction. By the reasoning in Part $(1)$, the total space of the pulled back bundle continues to be PSC. If $\dimmc(\widetilde{M}) = 3$, then the the macroscopic dimension of this cover is also $3$. If $M$ is enlargeable, then the cover is also enlargeable: it admits a degree $2$ map to $M$ and it is spin since all orientable $3$-manifolds are spin. Hence, the cover is also enlargeable. Therefore, we may also assume without loss of generality that the bundle $E \to M$ is fiberwise orientable.
    \end{enumerate}
   Let $P_1, \cdots, P_\ell$ be all the $3$-manifolds with finite fundamental group appearing in the Kneser-Milnor prime decomposition of $M$. By the Kurosh subgroup theorem, the kernel of the homomorphism $\varphi : \pi_1(M) \to \pi_1(P_1) * \cdots * \pi_1(P_\ell) \to \pi_1(P_1) \times \cdots \times \pi_1(P_\ell)$
   is torsion-free. Replacing $M$ with the finite cover corresponding to $\ker(\varphi)$ if necessary, we may assume without loss of generality that $\pi_1(M)$ is torsion-free. By Proposition \ref{prop-asphsum}, $M$ must contain at least one aspherical summand. Thus, $M = K \# N$ where $K$ is an orientable aspherical $3$-manifold. If the inclusion of the fiber $S^1 \hookrightarrow E$ is homotopically nontrivial, it follows from \cite[Theorem 1.9]{He25} that $E$ is not PSC. Therefore, let us assume the fibers of $E \to M$ are $\pi_1$-null in the total space of $E$. Thus, $\pi_1(E) \cong \pi_1(M)$ is a torsion-free $3$-manifold group. As such, $\pi_1(E)$ has cohomological dimension at most $3$. Furthermore, $\pi_1(E)$ satisfies the Strong Novikov Conjecture \cite[Theorem 1.1]{MOP08}\footnote{This theorem proves that the Baum-Connes Conjecture with coefficients holds for an oriented $3$-manifold group $G$. The discussion on \cite[pg.~5]{MOP08} shows that the Strong Novikov Conjecture for $G$ follows if the Baum-Connes Conjecture with coefficients holds for $G$.}. Next, observe that as $M$ is a closed orientable $3$-manifold, it is spin. Choosing a connection on $p : E \to M$, we obtain a decomposition $TE \cong V \oplus H$ into the vertical and horizontal subbundles. Note that $H \cong p^*TM$ and $V$ is the trivial line bundle over $E$. Therefore, $w_2(TE) = w_2(H \oplus V) = w_2(H) = p^* w_2(TM) = 0$. As $M$ is orientable and $p : E \to M$ is fiberwise orientable, the total space $E$ is also orientable. Since $w_2(TE) = 0$, this implies the total space $E$ is spin. It now follows from \cite[Main Theorem, Remark 1.7]{bol09} that if $E$ is PSC, the classifying map $f : E \to K(\pi_1(E), 1)$ is homotopic into the $2$-skeleton of $K(\pi_1(E), 1)$. However, this contradicts Lemma \ref{lem-nullfibers}. Thus, $E$ is not PSC. This concludes the proof.
\end{proof}

\section{Further examples and remarks}\label{sec-eg}

In this section, we provide several examples of manifolds satisfying the conclusions of Theorem \ref{main}. We begin with the simplest example given by a symplectic torus: 

\begin{exam}[Torus]\label{eg-torus}
    Let $M = T^{2n}$ for $n \geq 3$ and let $\omega$ be the standard symplectic form on $T^{2n} = \mathbf{R}^{2n}/\mathbf{Z}^{2n}$, normalized to be integral. Since $\pi_2(T^{2n}) = 0$, $(M, \omega)$ is a symplectically aspherical manifold. Furthermore, $\pi_1(T^{2n}) = \mathbf{Z}^{2n}$ is an amenable group. Moreover, $M = T^{2n-2} \times T^2$ as symplectic manifolds, and $T^{2n-2}$ is enlargeable. Therefore, $(M, \omega)$ satisfies the hypothesis of Parts $(3)$ and $(4)$ of Theorem \ref{main}. 
\end{exam}

We remark that using Example \ref{eg-torus}, we can produce examples of macroscopically large manifolds $Z$ admitting circle bundles $p : E \to Z$ such that the total space $E$ is macroscopically small, yet the asymptotic dimension of $\pi_1(E)$ is large. This indicates that the macroscopic smallness of $E$ cannot be deduced from just the fundamental group.\footnote{For any metric space $X$, $\dimmc(X) \leq \mathrm{asdim}(X)$. For more information on the asymptotic dimension, we refer the reader to the survey \cite{BeDr}.} 

\begin{exam} 
    Let $m \geq 2$, and $Z^{2m} \subset (T^{2n}, \omega)$ be a symplectic submanifold obtained by repeated applications of Theorem \ref{thm-dondiv}. Hence, by Proposition \ref{zpi12lemma}, $\pi_1(Z) \cong \ZZ^{2n}$. By Part $(2)$ of Theorem \ref{main}, there exists a circle bundle $S^1 \hookrightarrow E \stackrel{p}{\to} Z$ such that ${\dimmc(\widetilde{E}) \leq m+1}$. From the homotopy long exact sequence for fiber bundles, we obtain an exact sequence
    $${\pi_2(Z)} \to \pi_1(S^1) \cong \ZZ \to \pi_1(E) \to \pi_1(Z) \to 0$$
    {If $m \geq 3$, then $\pi_2(Z) \cong \pi_2(T^{2n}) \cong 0$ by Proposition \ref{zpi12lemma}. Therefore, $\pi_1(E)$ is an infinite cyclic extension of $\mathbf{Z}^{2n}$.} Therefore, the asymptotic dimension of $\pi_1(E)$ is
    $\mathrm{asdim}(\pi_1(E)) = 2n+1$ {by \cite[Theorem 66]{BeDr}, since $\pi_1(E)$ is polycyclic with Hirsch length $2n+1$}. Likewise, for the virtual cohomological dimension, we have $\mathrm{vcd}(\pi_1(E)) = 2n+1$. {If $m = 2$, then possibly $\pi_2(Z) \to \pi_1(S^1) \cong \mathbf{Z}$ has non-zero image. In that case, $\pi_1(E)$ is a finite cyclic extension of $\mathbf{Z}^{2n}$. Then $\mathrm{asdim}(\pi_1(E)) = 2n$ by \cite[Corollary 54]{BeDr} and since $\mathbf{Z}^{2n}$ is a finite-index subgroup of $\pi_1(E)$, $\mathrm{vcd}(\pi_1(E)) = \mathrm{vcd}(\mathbf{Z}^{2n}) = 2n$.} 
\end{exam}

Next, we give examples of symplectic manifolds which are enlargeable. These can be used as an input in Part $(4)$ of Theorem \ref{main}. To construct these examples, we use the recent work of Fine and Panov \cite{finepanov}, which shows that every even dimensional closed oriented manifold is dominated by a symplectic manifold.

\begin{exam}[Enlargeable symplectic manifolds]
    Let $N^{2n}$ be an even dimensional closed enlargeable manifold, for $n \geq 3$. By \cite[Theorem 1]{finepanov}, there exists a symplectic manifold $(M^{2n}, \omega)$ and a positive degree map $f : M \to N$. {We may arrange $c_1(M)$ to be an even multiple of $[\omega]$ by ensuring that in the proof of \cite[Theorem 1]{finepanov}, the chain of symplectic submanifolds of the twistor space of $N$, is selected so that all (resp. all but one) of them is a Donaldson divisor that is Poincar\'{e} dual to an even multiple of the symplectic form in the previous submanifold, depending on whether the parity of $n$ is even (respectively, odd). This ensures $c_1(M)$ is even by the remark at the end of the article \cite[pg.~3]{finepanov} regarding the first Chern class of the twistor space of $N$, together with the expression (\ref{eq-3}) in Proposition \ref{prop-unique} for the first Chern class of a Donaldson divisor of a symplectic manifold.} Therefore, we have
    $$w_2(M) \equiv c_1(M) \equiv 0 \pmod{2}$$
    Thus, $M$ is spin. Since $M$ admits a positive degree map to the enlargeable manifold $N$, $M$ must itself be enlargeable. Therefore, for any surface with an area form $(S_g, \omega_0)$, the symplectic manifold $(M \times S_g, \omega \oplus \omega_0)$ satisfies the hypothesis of Part $(4)$ of Theorem \ref{main}.
\end{exam}

We now give some examples which do not satisfy the hypotheses of Parts (3) and (4) in Theorem \ref{main}, nonetheless we observe that they satisfy the conclusions using similar methods.  

\begin{exam}[Product with a $2$-sphere]\label{eg-prodS2} {Let $n \geq 2$}. Let $(M^{2n}, \omega)$ be a closed integral symplectic manifold which is symplectically aspherical such that $\pi_1(M)$ is amenable. Let $\omega_0$ denote the standard area form on $S^2$. Then $(M \times S^2, \omega \oplus \omega_0)$ is not symplectically aspherical, hence it fails the hypothesis of Part $(3)$ of Theorem \ref{main}.

Nevertheless, let $Z^{2n} \subset (M \times S^2, \omega \oplus \omega_0)$ be a Donaldson divisor. Consider the map $f : M \times S^2 \to M$ given by the projection to the first factor. By Proposition \ref{zpullm}, $f|_Z : Z \to M$ is a positive degree map. Since $M$ is symplectically aspherical, it is rationally essential due to Proposition \ref{prop-sympasphess}. Since $f|_Z$ has positive degree, $Z$ must also be rationally essential. By Proposition \ref{zpi12lemma}, $\pi_1(Z) \cong \pi_1(M)$. Therefore, $Z$ is a rationally essential manifold with amenable fundamental group. {By \cite[Theorem 1.6]{dr11}, $\dimmc(\widetilde{Z}) = 2n$}.
\end{exam}

\begin{exam}[Symplectic blow-up]
    Let $M = T^{2n} \# \overline{\mathbf{CP}}^n$ for $n \geq 3$. By the symplectic blowup construction {\cite[pg 324-325]{msbook}}, $M$ admits a symplectic structure $\omega$. Let ${\overline{\mathbf{CP}}^{n-1} \subset M}$ denote the exceptional divisor. Then, $\overline{\mathbf{CP}}^1 \subset {\overline{\mathbf{CP}}^{n-1}}\subset (M, \omega)$ is a symplectic submanifold. Therefore,
    $$\int_{\overline{\mathbf{CP}}^1} \omega \neq 0$$
    Thus, $(M, \omega)$ is not symplectically aspherical. In addition, $M$ is not a product of an enlargeable manifold with a surface. Therefore, $(M, \omega)$ does not satisfy the hypothesis of either Part $(3)$ or Part $(4)$ in Theorem \ref{main}.
    
    Nevertheless, let $Z^{2n-2} \subset (M, \omega)$ be a Donaldson divisor with $[Z] = \mathrm{PD}(k [\omega])$, for some $k \gg 1$. Consider the map $f : M \to T^{2n-2}$ given by composing the blow-down map $q : M \to M/{\overline{\mathbf{CP}}^{n-1}} \cong T^{2n}$ with the projection $T^{2n} \cong T^{2n-2} \times T^2 \to T^{2n-2}$. For $z \in T^{2n-2}$, 
    $$(f|_Z)^{-1}(z) = f^{-1}(z) \cap Z = q^{-1}(\{z\} \times T^2) \cap Z$$
    Choose $z \in T^{2n-2}$ such that $\{z\} \times T^2 \subset T^{2n}$ is disjoint from a neighborhood of the blown up point $\{z_0\} = q(\overline{\mathbf{CP}}^n) \subset T^{2n}$. By the construction of symplectic blow-ups, $q$ restricts to a symplectomorphism on the complement of a neighborhood of the exceptional divisor ${\overline{\mathbf{CP}}^{n-1}} \subset M$. Therefore, $T^2 \cong q^{-1}(\{z\} \times T^2) \subset (M, \omega)$ is a symplectic submanifold intersecting $Z^{2n-2}$ transversely. Therefore, the oriented intersection number is given by:
    $$\# (q^{-1}(\{z\} \times T^2) \cap Z) = k \int_{q^{-1}(\{z\} \times T^2)}  \omega > 0$$
    Thus, $f|_Z : Z \to T^{2n-2}$ is a map of positive degree. Hence, $Z$ is rationally essential. By Proposition \ref{zpi12lemma}, $\pi_1(Z) \cong \pi_1(M) \cong \mathbf{Z}^{2n}$ which is amenable. Therefore, by {\cite[Theorem 1.6]{dr11}}, $\dimmc(\widetilde{Z}) = 2n-2$. Additionally, if $n = 3$ (resp. $n = 4$), then $Z$ admits a positive degree map to $T^4$ (resp. $T^6$). In that case, $Z$ is not PSC \cite{sydescent}.  
\end{exam}

In dimension $4$, it is possible to produce many more examples of non-PSC manifolds admitting circle bundles with total space PSC using a deep result of Taubes \cite{taubes}. 

\begin{exam}
    Let $(M^4, \omega)$ be a closed integral symplectic $4$-manifold with $b_2^+(M) > 1$, and let $\omega_0$ denote an integral area form on the surface $S_g$. Let $Z^4 \subset (M \times S_g, \omega \oplus \omega_0)$ be a Donaldson divisor. By Proposition \ref{zpullm}, the projection $f : Z \to M$ has positive degree. Therefore, the induced map 
    $$f^* : H^*(M; \QQ) \to H^*(Z; \QQ)$$
    is an injection of quadratic spaces, where both the domain and range are equipped with the quadratic form given by cup product. Therefore, since $b_2^+(M) > 1$, we must have $b_2^+(Z) > 1$ as well. Since $Z$ is a symplectic $4$-manifold with $b_2^+(Z) > 1$, by Taubes' theorem \cite{taubes}, it follows that $Z$ has non-vanishing Seiberg-Witten invariants. Therefore, $Z$ is not PSC. Nevertheless, by Part $(1)$ of Theorem \ref{main}, $Z$ admits a circle bundle with total space being PSC.
\end{exam}

We give an example which illustrates that macroscopic dimensions of manifolds having diffeomorphic universal covers need not be the same. For more surprising examples, we refer the reader to Kastenholz-Reinhold \cite{kr}.

\begin{exam}
    Let $M = T^{2n}$ for $n \geq 3$ and $\omega$ be the standard symplectic form on the torus, normalized to be integral. Let $Z \subset (M, \omega)$ be a Donaldson divisor with $[Z] = \mathrm{PD}(k [\omega])$. By Part $(2)$ of Theorem \ref{main}, there exists a circle bundle $S^1 \hookrightarrow E \stackrel{p}{\to} Z$ such that $\dimmc(\widetilde{E}) \leq n$. 
    Let $\widehat{E} \to \widetilde{Z}$ denote the circle bundle obtained by pulling back the bundle $p : E \to Z$ along the universal cover $q : \widetilde{Z} \to Z$. The resulting bundle is classified by an element in
    $$H^2(\widetilde{Z}; \ZZ) \cong \mathrm{Hom}(H_2(\widetilde{Z}; \ZZ), \ZZ) \cong \mathrm{Hom}(\pi_2(\zcov), \ZZ) \cong \mathrm{Hom}(\pi_2(Z), \ZZ)$$
    More precisely, it is classified by the cohomology class of the $2$-form ${k \cdot q^*\omega}$. However, since $(Z, \omega)$ is symplectically aspherical by Example \ref{eg-sympash}, for any map $f : S^2 \to \widetilde{Z}$, 
    $$\int_{S^2} f^* (q^*\omega) = \int_{S^2} (q \circ f)^*\omega = 0$$
    Therefore, $[q^* \omega] = 0 \in H^2(\widetilde{Z}; \ZZ) \cong \mathrm{Hom}(\pi_2(Z), \ZZ)$. This implies that the bundle $\widehat{E} \to \widetilde{Z}$ is trivial, i.e., $\widehat{E} \cong \widetilde{Z} \times S^1$. Since the total space of $\widehat{E}$ is a cover of the total space of $E$, we obtain that $\widetilde{Z} \times S^1$ covers the total space of $E$. In particular, $\widetilde{E}$ is diffeomorphic to $\widetilde{Z} \times \RR$. Note that $\dimmc(\widetilde{Z}) = 2n-2$, by Part $(3)$ of Theorem \ref{main}. Thus, $E$ and $Z \times S^1$ have diffeomorphic universal covers, but nevertheless have different macroscopic dimensions.
\end{exam}

Recall the discrepancy (see, \cite{drmoscow}) between the two notions of macroscopic dimension introduced in Definition \ref{def-md}. We close this section by giving examples of manifolds $Z$ with $\dim_{\mathrm{mc}}(\zcov) 
= \dim(\zcov)$ (as opposed to $\dimmc(\zcov) = \dim(\zcov)$, as discussed throughout the article), such that there exists macroscopically small and PSC circle bundles over $Z$.

\begin{exam}[md-large divisors] \label{examplemdlarge}
Let $(T^{2n}, \omega)$, ${n \geq 2}$ be a symplectic torus as in Example \ref{eg-torus}. Consider the product $(T^{2n} \times S^2, \omega \oplus \omega_0)$, as in Example \ref{eg-prodS2}. Let $Z^{2n} \subset (T^{2n} \times S^2, \omega \oplus \omega_0)$ be a Donaldson divisor. Then there is a positive degree map $f : Z \to T^{2n}$ which induces an isomorphism on $\pi_1$, by Example \ref{eg-prodS2}. By lifting $f$ to the universal cover, we obtain a map with positive local degree
$$\widetilde{f} : \widetilde{Z} \to \widetilde{T^{2n}} \cong \mathbf{R}^{2n}$$
Therefore, $\widetilde{f}$ induces a non-zero homomorphism between the locally finite homology groups
$\widetilde{f}_* : H_{2n}^{lf}(\widetilde{Z}; \ZZ) \to H_{2n}^{lf}(\RR^{2n}; \ZZ)$. Consequently, $\widetilde{f}_*[\widetilde{Z}] \neq 0$. Since $\pi_1(Z) \cong \ZZ^{2n}$, we have $T^{2n} = K(\pi_1(Z), 1)$. Therefore, by \cite[Theorem 5.4]{drmoscow}, $\dim_\mathrm{mc}(\widetilde{Z}) = 2n$. By Parts $(1)$ and $(2)$ of Theorem \ref{main}, we obtain circle bundles $E \to Z$ such that $\dim_{\mathrm{mc}}(\widetilde{E}) \leq \dimmc(\widetilde{E}) \leq n + 1$ and $E$ is PSC.
\end{exam}

\bibliographystyle{alpha}
\bibliography{pscbundles-bib}
\end{document}